\RequirePackage{fix-cm}
\documentclass[english]{scrartcl}
\usepackage[T1]{fontenc}
\usepackage[latin9]{inputenc}
\usepackage{geometry}
\usepackage{amsthm}
\usepackage{amsmath}
\usepackage{amssymb}

\makeatletter
\numberwithin{equation}{section}
\numberwithin{figure}{section}
\@ifundefined{lettrine}{\usepackage{lettrine}}{}

\usepackage{ifxetex}\ifxetex
\usepackage{fontspec}
\else
\fi

\usepackage{fixltx2e}\usepackage{fancyhdr}

\makeatother

\usepackage{babel}
\begin{document}

\title{Consequences of a Gödel's misjudgment}

\author{Giuseppe Raguní}

\date{UCAM - Universidad Católica de Murcia, Spain - graguni@ucam.edu}
\maketitle
\begin{abstract}
The fundamental aim of the paper is to correct an harmful way to interpret
a Gödel's erroneous remark at the Congress of Königsberg in 1930.
Despite the Gödel's fault is rather venial, its misreading has produced
and continues to produce dangerous fruits, as to apply the incompleteness
Theorems to the \emph{full} second-order Arithmetic and to deduce
the semantic incompleteness of its language by these same Theorems.
The first three paragraphs are introductory and serve to define the
languages \emph{inherently semantic} and its properties, to discuss
the consequences of the expression order used in a language and some
question about the semantic completeness: in particular is highlighted
the fact that a non-formal theory may be semantically complete despite
using a \emph{language} semantically incomplete. Finally, an alternative
interpretation of the Gödel's unfortunate comment is proposed.

\end{abstract}

\section{Formal systems}

Often the adjective \emph{formal} is abused in violation of the original
meaning that can rightly be called \emph{Hilbertian}. Actually, it
is common understanding ``formal system'' as a synonym for ``axiomatic
system'', although no one doubts the need for non-formal axioms in
Mathematics. This superficiality, together with the wrong meaning
given to certain affirmations of Gödel, produces the serious mistakes
that we will show. 

The best definition of a formal system is probably that one given
by Lewis in 1918%
\footnote{C. I. Lewis, \emph{A Survey of Formal Logic}, Berkeley University
of California Press (1918), p. 355. %
}:
\begin{quotation}
A mathematical system is any set of strings of recognizable marks
in which some of the strings are taken initially and the remainder
derived from these by operations performed according to rules which
are independent of any meaning assigned to the marks.
\end{quotation}
So, a framework in which, starting by certain strings of characters
with no meaning (the formal axioms), are produced (deducted) other
strings with no meaning (the theorems) by use of operations (deductive
rules) which never make use of any meaning assigned to any mark. To
get from that sort of ``abacus'' a real scientific discipline, we
must \emph{interpret} - if this is possible - the characters and strings
so that: a) the axioms are \emph{true}; b) the deductive rules are
\emph{sound}, ie capable of producing only \emph{true} meaningful
strings, starting by\emph{ true} meaningful strings. Such an interpretation
is said \emph{sound} \emph{model}. This conceptualization, revolutionary
for the times, is also present in the contemporary Bernays and Post%
\footnote{See, for example, the Dreben and Heijenoort introduction to first
Gödel's works in \emph{Kurt Gödel collected works}, ed. Feferman \emph{et
al.}, Oxford univ. press (1986), vol. 1, pp. 44-48.%
} as well - perhaps for the very first time - in Hilbert, as it can
be noted in many passages of his correspondence with Frege (1899-1900):
\begin{quotation}
{[}...{]} It is surely obvious that every theory is only a scaffolding
or schema of concepts together with their necessary relations to one
another, and that the basic elements can be thought of in any way
one likes. If in speaking of my points I think of some system of things,
e.g. the system: love, law, chimney-sweep\dots{} and then assume all
my axioms as relations between these things, then my propositions,
e.g. Pythagoras' theorem, are also valid for these things. In other
words: any theory can always be applied to infinitely many systems
of basic elements.

{[}...{]} a concept can be fixed logically only by its relations to
other concepts. These relations, formulated in certain statements
I call axioms, thus arriving at the view that axioms\dots{} are the
definitions of the concepts. I did not think up this view because
I had nothing better to do, but I found myself forced into it by the
requirements of strictness in logical inference and in the logical
construction of a theory%
\footnote{From two Hilbert's letters to Frege, in: G. Frege, \emph{Philosophical
and Mathematical Correspondence}, ed. G. Gabriel, \emph{et al.} Oxford:
Blackwell Publishers (1980).%
}.
\end{quotation}
Incidentally, the fact that none of the mentioned authors uses rigorously
the adjective \emph{formal} is not surprising: since for all of them
the forthright belief (or hope) was that \emph{every} axiomatic theory
was formalizable, the formal dressing was nothing more than the \emph{correct
way} of presenting any axiomatic discipline. Precisely Gödel is one
of the first to use the term carefully and accurately for his attention
to distinguish meticulously between Mathematics and Metamathematics
and his consciousness of expressive limits for the formality (which
he self will help to shape). As a matter of fact, even without resorting
to the feasibility of the condition of categoricity (which only applies
in non-formal axiomatic systems%
\footnote{If we limit ourselves to the consideration of theories with models
of infinite cardinality.%
}), the non-formal axiomatic theories are necessary to express, and
possibly decide, all that the Mathematics was thought to deal. 

The concept of formality requires at least a clarification. The question
is simply whether the \emph{formal} systems coincide with \emph{recursively}
(or \emph{effectively}, assuming the Church-Turing Thesis) \emph{axiomatizable}
systems. When the works of Church and, especially, Turing%
\footnote{A. M. Turing, \emph{On computable numbers, with an application to
the Entscheidungsproblem}, Proceedings of the London Mathematical
Society (2) 42 (1937), pp. 230-265.%
} convinced Gödel that the machine models considered by them were,
in fact, completely equivalent to his general recursive functions
(ie convinced him of the Church-Turing Thesis), he, in a footnote
added in 1963 to the text of his most famous work, wrote:
\begin{quotation}
In consequence of later advances, in particular of the fact that due
to A. M. Turing's work a precise and unquestionably adequate definition
of the general notion of formal system can now be given, a completely
general version of Theorems VI and XI {[}\emph{the two incompleteness
Theorems}{]} is now possible. That is, it can be proved rigorously
that in \emph{every} consistent formal system that contains a certain
amount of finitary number theory there exist undecidable arithmetic
propositions and that, moreover, the consistency of any such system
cannot be proved in the system%
\footnote{K. Gödel\emph{, On formally undecidable propositions of Principia
mathematica and related system I, }in\emph{ Kurt Gödel collected works},
\emph{op. cit.}, vol. 1 p. 195.%
}.
\end{quotation}
clarifying, in a footnote of this same text, that:
\begin{quotation}
In my opinion the term ``formal system'' or ``formalism'' should
never be used for anything but this notion. {[}...{]}
\end{quotation}
The Gödel's proposal was therefore to coincide, by definition, the
\emph{effectively axiomatizable} systems with the \emph{formal} systems.
Unfortunately, this suggestion was not successful. Certainly, it is
just a matter of \emph{convention}; but if we maintain the original
definition of \emph{formal}, the two concepts are different. Consider
a theory that deduces only by the classical deductive rules. Since
these rules are mechanizable, the only obstacle for a machine, in
simulating the system, is the specification of the set of axioms:
if it is mechanically reproducible, the system will be \emph{effectively
axiomatizable}. On the other hand, on the basis of definition of \emph{formality},
the chains to call \emph{axioms} undoubtedly must be exhibitable but
can be specified in arbitrary way. So, \emph{effective axiomatizability}
implies \emph{formality} but not vice versa.

An important example of this case can be built by the formal Peano
Arithmetic (\emph{PA}). Assuming that such theory is consistent, we
can add as new axioms the class of its \emph{true} statements in the
intuitive (or \emph{standard}) model. By construction, it forms a
syntactically complete system (\emph{PAT}) which, as a result of the
first incompleteness Theorem, cannot be effectively axiomatizable.
Nevertheless this system can still be called \emph{formal}: it is
necessary to give meaning to the sentences just to determine whether
or not they are axioms. Once done, every statement can be reconverted
to a meaningless string, since the system can deduce its theorems
without making use of significance. So, really exists a formal system
able to solve the halting problem: cold comfort.

\section{Semantic languages and non-formal systems}

Consider an arbitrary language that, as normally, makes use of a countable%
\footnote{Finite, as the usual alpha-numeric symbols, or, to generalize, infinite-countable.
In this paper we will use either ``countable'' or ``enumerable''
with the same meaning, ie to indicate that there exists a biunivocal
correspondence between the considered set and the set of the natural
numbers. %
} number of characters. Combining these characters in certain ways,
are formed some fundamental strings that we call \emph{terms} of the
language: those collected in a dictionary. When the terms are semantically
interpreted, ie a certain meaning is assigned to them, we have their
distinction in adjectives, nouns, verbs, etc. Then, a proper grammar
establishes the rules of formation of sentences. While the terms are
finite, the combinations of grammatically allowed terms form an infinite-countable
amount of possible sentences.

In a non-trivial language, the meaning associated to each term, and
thus to each expression that contains it, is not always unique. The
same sentence can enunciate different things, so representing different
\emph{propositions}. For example, the same sentence ``it is a plain
sailing'' has a different meaning depending on the circumstances:
at board of a ship or in the various cases with figurative sense.
How many meanings can be associated to the same term? That is: how
many different propositions, in general, can we get by a single sentence?
The answer, for a normal semantic language, may be amazing.

Suppose we assign to each term a finite number of well-defined meanings.
We could then instruct a computer to consider all the possibilities
of interpretation of each term. The computer, to simplify, may assign
all the different meanings to an equal number of distinct new terms
that it has previously defined. For example, it might define the term
``f-sailing'' for the figurative use of ``sailing'' (supposed
unique). The machine would then be able, using the grammar rules,
to generate all the infinite-countable propositions. In this case
we will say that, in the specific language, the meaning has been \emph{deleted}%
\footnote{Just a concise choice rather than the more correct \emph{mechanically
reproduced}.%
}\emph{.} More generally we have this case when the different meanings
allowed for each term are \emph{effectively enumerable}: even in the
case of an infinite-countable amount of meanings, the computer can
define an infinite-countable number of new terms and associate only
one meaning to each term in order to establish a biunivocal correspondence
between sentences and propositions. So, the machine could list all
them by combination.

Hence, by definition, we will say that a language is \emph{inherently
semantic} (ie with a non-eliminable meaning) if it uses at least one
term with an amount not effectively enumerable of meanings; with the
possibility, which we will comment soon, that this quantity is even
uncountable. From the fact that a sentence represents more than one
proposition if and only if it contains at least one term differently
interpreted, it follows an equivalent condition for the inherent semanticity:
a language is inherently semantic if and only if the set of all possible
propositions is not effectively enumerable. 

A first important example of such a language was considered in the
previous section. The fact that the axioms of \emph{PAT} are not effectively
enumerable proves that in the expression ``true statement in the
standard model'' the term \emph{true} has got an amount not effectively
enumerable (although enumerable) of distinct meanings. So, the phrase
belongs to an inherently semantic language.

When the set of all possible propositions is enumerable but not effectively
enumerable, still it is possible to define an amount infinite-countable
of new terms and to associate only one meaning to each term (so re-establishing
a biunivocal correspondence between sentences and propositions); but
this operation cannot be performed by a machine. Go back to the example\emph{
}of \emph{PAT}: an inherently semantic phrase is used to define all
its new axioms, but then every axiom is formulated by its unique symbolic
representation in \emph{PA}. Every proposition of \emph{PAT} corresponding
to the expression ``true statement in the standard model'' is explicitly
replaced with an appropriate formula of \emph{PA}. So, in this case,
the ``definition of new terms'' consists precisely in using the
formulas of \emph{PA} to express the propositions of \emph{PAT}. 

Now consider the case of an inherently semantic language in which
all possible propositions are not enumerable, ie in which there exists
at least one term with an uncountable quantity of meanings. Due to
the uncountability of properties of the standard natural numbers (\emph{N}),
ie of the set of all subsets of \emph{N}, \emph{P(N)}, this case is
especially interesting, because only a language of this kind is capable
of expressing these properties. This time it is not possible to define
new terms so to associate to them an unique meaning, because the number
of all possible strings is only countable%
\footnote{Sometimes the uncountability is interpreted simply as a ``non-formalizable
countability'', ie achieved by an inherently \emph{metamathematical}
connection. Even according to this view, however, this ``countability''
cannot be capable to assign a different symbolic string to each element
of the uncountable set: otherwise, nothing could prohibit to formalize
this function inside the formal Set Theory.%
}. Not even we can reaccommodate the things so to re-associate a quantity
at most infinite-countable of interpretations to each term, because
in this way we still could get only a \emph{countable} total amount
of sentences. \emph{There is no way to avoid that at least one term
conserves an uncountable number of interpretations}.

Although that is what really happens in every usual natural language,
this feature, at first, might surprise or even be considered unacceptable.
Undoubtedly, all the meanings that \emph{ever will be assigned} to
any settled word are only a countable number: indeed, finite! But
these meaning cannot be specified once and for all. The fact remains
that the \emph{possible} interpretations of the term vary inside an
infinite collection; moreover, a collection not limited by any prefixed
cardinality. Some classic paradoxes can be interpreted as a confirmation
of this property. The Richard's one%
\footnote{Where firstly it is admitted that all possible semantic definitions
of the real numbers stay in a countable array and, then, one can define,
by a diagonal criterion, a real number that is not present in the
array.%
}, for example, can be interpreted as a meta-proof that the semantic
definitions are not countable, ie that they are \emph{conceivably}
able to define each element of a set with cardinality greater than
the enumerable one (and therefore each real number). The proper technique
used in \emph{diagonal} \emph{reasoning}, reveals that the natural
language is able to adopt different semantic levels (or contexts)
``looking from the outside'' what was ``before'' defined, namely,
what was previously said by the same language. Identical words used
in different contexts have a different meaning and for the number
of contexts, including nested, there is no limit. On the other hand,
the Berry's paradox clearly shows that a finite amount of symbols,
differently interpreted, is able to define an infinite amount of objects.
Here the key of the argument is again the use of two different contexts
to interpret the verb ``define''.

Finally, we consider a consistent arbitrary axiomatic system (\emph{S\textsubscript{A}}).
We wish that \emph{S\textsubscript{A}} is able to express and possibly
decide%
\footnote{By ``decide a sentence'', we mean to conclude that it or its negation
is a theorem.%
} all the properties of the natural numbers. In particular it must
be capable to distinguish the properties one from the other, in order
to deduce, in general, different theorems starting by different properties.
Admitting, as normally, that \emph{S\textsubscript{A}} makes use
of a countable number of symbols, if we interpret the theory in a
certain conventional model, this will associate a single meaning to
each sentence. So, the interpreted sentences will be only an enumerable
amount and they cannot express all the properties of the natural numbers.
Therefore, \emph{the only possibility} \emph{is to consider a non-conventional
model} of the system, able to assign more than a single meaning (indeed
a quantity at least 2\textsuperscript{$\aleph_{0}$}) to at least
one sentence; and, then, able to verify the remaining requirements
for a normal model%
\footnote{As a result, this model cannot be confused with a conventional uncountable
model of a formal, and maintained formal, system. This one, considered
for example for \emph{PA}, will continue to assign \emph{only one}
meaning to each sentence of the system: an uncountable amount of exceeding
elements of the universe will have no representation in \emph{PA}.%
}. That is what really happens in the so-called \emph{full} second
order Arithmetic (\emph{FSOA}): its standard model, which can be proved
unique under isomorphism, is of this kind. As we said, this axiomatic
system\emph{ is} \emph{non-formal}. Definitively, \emph{no kind of
interpretation can allow an axiomatic system to express and study
all the properties of the natural numbers in compliance with formality}.

\section{Questions about the expression orders}

The first-order predicate calculus of Logic%
\footnote{Defined by Russell and Whitehead in \emph{Principia Mathematica},
Cambridge University Press (1913).%
}, or briefly \emph{first-order classical Logic}, is a formal system
which constitutes the structural core (ie the\emph{ language}) of
the ordinary axiomatic theories. In it, the existential quantifiers
$\forall$ and $\exists$%
\footnote{Really, only one of the two is strictly necessary.%
} only can range over variable-elements of the universe (\emph{U})
of the model. In 1929 Gödel proved the semantic completeness of this
theory, namely, that all the \emph{valid} (ie true in every model)
sentences are theorems%
\footnote{K. Gödel, \textit{op. cit.,} doc\emph{. 1929 }and \emph{1930}, vol.
1, p. 61 and p. 103.%
}. The widest generalization of the Gödel's result is due to Henkin%
\footnote{L. Henkin, \emph{The completeness of the first-order functional calculus},
The journal of symbolic logic (1949), n.14, pp. 159-166.%
} and extends the semantic completeness to \emph{every formal system}.

The second-order predicate calculus of Logic, or \emph{second-order
classical Logic}, extends the use of existential quantifiers to predicates,
ie to properties of the elements of \emph{U}. From the standpoint
of the Set theory, the predicate-variables vary inside the set of
all subsets of \emph{U}, that is, \emph{P(U)}. In non-trivial cases,
\emph{U} is infinite and consequently \emph{P(U)} always has got an
uncountable cardinality. If we consider an arbitrary non-trivial axiomatic
theory based on second-order classical Logic, we have thus three fundamental
cases:
\begin{enumerate}
\item The axioms and/or rules of the theory nevertheless ensure formality.
For this purpose it is necessary (although not sufficient) that they
limit the variability of the predicates within a \emph{countable}
subset of \emph{P(U)}; normally, this is achieved by means of appropriate
\emph{comprehension axioms}. This case is known as \emph{general,
or Henkin's, Semantics}.
\item There is no limit for variability of the predicates in\emph{ P(U)}.
This \emph{full} understanding is known as \emph{standard Semantics}.
In any case a non-formal system is obtained, since it is possible
to express an uncountable quantity of sentences.
\item The variability of predicates in \emph{P(U)} is restricted but not
enough to verify the formality. The sentences may be uncountable or
countable but, in the latter case, not \emph{effectively}.
\end{enumerate}
In the case 1) the theory always is semantically complete; what happens
in the remaining cases?

The \emph{FSOA} Arithmetic is right an important example of case 2).
In it, an induction principle, valid for ``any property of the natural
numbers'', is defined as axiomatic scheme. Definitely, this is a
\emph{full} understanding of the inductive rule, that implies uncountability
for the sentences of the theory. By the way, it is well known that
such \emph{full} understanding \emph{is necessary to achieve the categoricit}y%
\footnote{A full explanation, for example, in S. Shapiro, \emph{Foundations
without foundationalism: a case for second-order Logic}, Oxford univ.
press (2002), p. 82 and p. 112. %
}. By its categoricity, is immediate to show that\emph{ }the compactness
Theorem, which is a corollary of the semantic completeness Theorem,
does not hold for \emph{FSOA}. This allows the conclusion that \emph{the
language} of \emph{FSOA}, ie the \emph{full} second-order Logic, is
not semantically complete. A more general alternative for this conclusion
is obtained applying the Löwenheim-Skolem Theorem (\emph{L-S}). In
a simplified form that includes both versions (\emph{up} and \emph{down}),
the \emph{L-S} Theorem may be stated in this way: every system equipped
with a semantically complete language and with at least one infinite
model, admits models of any infinite cardinality. Since models of
different cardinality cannot be isomorphic, we have that for a system
with infinite models, \emph{categoricity and semantic completeness
of its language} \emph{are two complementary properties}, impossible
to be both satisfied.

At this point an explanation is necessary, which, surprisingly, I
never have found in any publication. Either by the invalidity of compactness
Theorem or by the \emph{L-S} Theorem, what one really concludes is
precisely the semantic incompleteness \emph{of the language} of \emph{FSOA},
ie of the \emph{full} second-order Logic;\emph{ and not of the system.}
\emph{Not every system that is expressed by the full second-order
Logic (or, more generally, by any other semantically incomplete language)
has necessarily to be semantically incomplete}%
\footnote{While, conversely, a system that uses a semantically complete language
always is semantically complete.%
}. In particular, for a categorical (non-trivial) system (\emph{which
always, by the L-S Theorem, uses a semantically incomplete} \emph{language}),
there is nothing to prohibit \emph{a priori} that it can be syntactically
(or semantically%
\footnote{For a categorical system it is immediate to prove that the semantic
and syntactic completeness are equivalent.%
}) complete. Indeed, such a system, which necessarily will be non-formal,
could deduce all the true sentences using an appropriate kind of inherently
semantic deduction%
\footnote{More in G. Raguní, \emph{I confini logici della matematica}, ed. Aracne,
Roma (2010), pp. 144-145 and 232-233.%
}.

On the other hand, a source of mistakes is undoubtedly the current
widespread tendency to categorize the axiomatic Theories looking at
the expression order (\emph{first order}, \emph{second order}, etc.)
without making the necessary distinctions. It is in line with the
disuse of the term \emph{formal} in its pure Hilbertian sense (so
used also by Gödel). Regardless to clarify if the semantics is \emph{full}
(or \emph{standard}) or \emph{general} (ie accomplishing the formality)
or intermediate (previous case 3), the systems of second (or more)
order are normally considered as those for which the semantic completeness,
and properties related to it, does not apply, as opposed to those
of first order. But, firstly, a first order classic system can own
\emph{proper} axioms that violate either the semantic completeness
or the same formality%
\footnote{See, for examples: M. Rossberg, \textit{First-Order Logic, Second-Order
Logic and Completeness}, Hendricks \emph{et al.} (eds.) Logos Verlag
Berlin (2004), \textit{\emph{on }}\textit{WEB}:

http://www.st-andrews.ac.uk/\textasciitilde{}mr30/papers/RossbergCompleteness.pdf%
}. Moreover, Henkin has shown that under conditions of formality, the
semantic completeness applies whatever the expression order. As we
noted in the previous section, the crucial aspect lies in the semantic
consequences of the cardinality of language of the theory, regardless
of the expression order. Now, the only peculiarity related to use
of the quantifiers $\exists$ and $\forall$ on predicative variables
(or super-predicative, etc.) is that when infinite models exist (so,
in non-trivial cases), these variable - if appropriate comprehension
axioms are missing - vary inside uncountable sets, so creating \emph{inherently
semantic} systems. But comprehension axioms can radically change the
things. A superficial interpretation of the Lindström's Theorem aggravates
the situation. It states that every classical theory expressed in
a semantically complete language \emph{can} be expressed in a first-order
language. This \emph{can} is not a \emph{must}. The theorem does not
prohibit semantic completeness or even formality for a system expressed
in an higher-order language%
\footnote{For the second order, you can have a semantically complete language
without it being formal, in one of the cases incorporated in previous
point 3.%
}. It only states that, when you have this case, the theory can be
re-expressed in a - simpler - first-order language. Certainly, this
property stands out the particular importance of the first-order language%
\footnote{A near property, on the other hand, is already evident thanks to the
expressive capacity of the formal Set Theory: namely, every \emph{formal}
system, since is fully representable in this theory - which is of
the first order - can be expressed at the first order.%
}. But this should not be radicalized. The grouping of the axiomatic
theories based on the expression order is, in general, misleading
their basic logical properties, unless the actual effect of the axioms
is considered. Because these properties are only a consequence of
the premises. The essential tool for classification remains the accomplishment
of the Hilbertian formality.

Starting from the \emph{FSOA}, the formality can be restored by limiting
the induction principle to the formally expressible properties, so
re-establishing an injective correspondence between formulas and properties:
in this way the system \emph{PA} is obtained, which however is unable
to express, and therefore also decide, infinite (namely, again 2\textsuperscript{$\aleph_{0}$},
since this number remains unchanged for subtraction of $\aleph_{0}$)
properties of the natural numbers. Now the induction principle is
not longer able to reject all the models not isomorphic to \emph{N},
so categoricity cannot be achieved%
\footnote{See, for example: S. Shapiro, \emph{op. cit.}, p. 112.%
}. Indeed, it can be shown in at least four different ways that \emph{PA}
has got non-standard models: 1. by the first incompleteness Theorem
together with the semantic completeness Theorem; 2. using the compactness
Theorem; 3. by the \emph{upward} \emph{L-S} Theorem; 4. as a consequence
of a Skolem's theorem of 1933 that concludes the non-categoricity
for \emph{PAT}%
\footnote{T. Skolem, \emph{Über die Unmoglichkeit...}, Norsk matematisk forenings
skrifter, series 2, n. 10, 73-82, reprinted in \emph{Selected works
in logic}, edited by Jens E. Fenstad, Univ. di Oslo (1970), pp. 345-354.%
}. 

Finally we recall that the incompleteness Theorems are valid for any
effectively axiomatizable (and therefore formal) system in which the
general recursive functions are definable. Therefore, they can be
applied to \emph{PA} but not to \emph{FSOA}. For \emph{PA}, the first
Theorem reveals a further limitation: if one admits its consistency,
even between the properties of the natural numbers that this theory
is able to express, there are infinite ($\aleph_{0}$) that it cannot
decide.

\section{From error to horror}

In the third volume of the aforementioned \emph{Kurt Gödel's collected
works}, published in 1995, are collected the unpublished writings
of the great Austrian logician. According to the editors, the document
\emph{{*}1930c} is, in all probability, the text presented by Gödel
at the Königsberg congress on September 6, 1930%
\footnote{Indeed, in the \emph{textual notes} of the volume is written (p. 439):
\begin{quotation}
The copy-text for \emph{{*}1930c} {[}...{]} was one of several items
in an envelope that Gödel labelled \textquotedblleft Manuskripte Korrekt
der 3 Arbeiten in Mo{[}nats{]}H{[}efte{]} + \emph{Wiener Vorträge
über die ersten zwei}\textquotedblright{} (manuscripts, proofs for
the three papers in \emph{Monatshefte} {[}\emph{1930}, \emph{1931},
and \emph{1933i}{]} plus Vienna lectures on the first two.) On the
basis of that label, \emph{{*}1930c} ought to be the text of Gödel\textquoteright s
presentation to Menger\textquoteright s colloquium on 14 May 1930
- the only occasion, aside from the meeting in Königsberg, on which
Gödel is known to have lectured on his dissertation results {[}...{]}.
Internal evidence, however, especially the reference on the last page
to the incompleteness discovery, suggests that the text must be that
of the later talk. Since no other lecture text on this topic has been
found, it may well be that Gödel used the same basic text on both
occasions, with a few later additions.\end{quotation}
}. 

In the first part of the document, Gödel presents his semantic completeness
Theorem, extended to the ``restricted functional calculus'' (with
no doubt%
\footnote{To be convinced, just consult the note 3 in K. Gödel, \textit{op.
cit.,} doc. \emph{1930}, vol. 1, p. 103.%
} identifiable with the first-order classical Logic), which an year
before he proved in his doctoral thesis. After this exposure, he adds:
\begin{quotation}
{[}...{]} If the completeness theorem could also be proved for the
higher parts of logic (the extended functional calculus), then it
could be shown in complete generality that syntactical completeness
follows from monomorphicity {[}categoricity{]}; and since we know,
for example, that the Peano axiom system is monomorphic {[}categorical{]},
from that the solvability of every problem of arithmetic and analysis
expressible in Principia mathematica would follow.

Such an extension of the completeness theorem is, however, impossible,
as I have recently proved {[}...{]}. This fact can also be expressed
thus: The Peano axiom system, with the logic of Principia mathematica
added as superstructure, is not syntactically complete%
\footnote{K. Gödel, \textit{op. cit.,} doc. \emph{{*}1930c}, vol. 3, pp. 27-29.%
}.
\end{quotation}
In summary, Gödel affirms that is impossible to generalize the semantic
completeness Theorem to the ``extended functional calculus''. In
fact, in this case also the Peano axiomatic system, structured with
the logic of \emph{Principia Mathematica} (\emph{PM}), would be semantically
complete. But since this theory is categorical, would follow that
it is also syntactically complete. But just this last thing is false,
as he - surprise - announces to have proved.

Now, regardless of what Gödel meant by ``extended functional calculus'',
\emph{this affirmation contains an error}. We have in fact two cases:
\begin{description}
\item [{a.}] If Gödel understands by ``Peano axiomatic system structured
with the logic of the \emph{PM}'' the formal theory \emph{PA} (\emph{or
any other formal arithmetical system}) the error is precisely to regard
it as categorical.
\item [{b.}] If, however, he alludes to the only categorical arithmetic,
that is \emph{FSOA}, then Gödel errs applying to it his first incompleteness
Theorem. 
\end{description}
Of the two, just the second belief has been consolidating but - shockingly
- without reporting the error. Rather, exalting the merit of having
detected for the first time the semantic incompleteness of the \emph{full}
second-order Logic. Surely, to forming this opinion has been important
the influence of the following sentence contained in the second edition
(1938) of \emph{Grundzüge theoretischen der Logik} by Hilbert and
Ackermann%
\footnote{English translation by L. M. Hammond \emph{et al.} in \emph{Principles
of Mathematical Logic}, Chelsea, New York (1950), p. 130.%
}:
\begin{quotation}
Let us remark at once that \emph{a complete axiom system for the universally
valid formulas of the predicate calculus of second order does not
exist.} Rather, as K. Gödel has shown {[}K. Gödel, \emph{Über formal
unentscheidbare Sätze der Principia Mathematica und verwandter systeme},
Mh. Math. Physik Vol. 38 (1931){]}, for any system of primitive formulas
and rules of inference we can find universally valid formulas which
cannot be deduced.
\end{quotation}
The echo of the Gödel's unfortunate words at the Congress pushes the
authors (probably Ackermann, given the age of Hilbert) to attest that
the first incompleteness Theorem concludes precisely the semantic
incompleteness of the second-order Logic! False. And furthermore,
this conclusion cannot be derived by that Theorem. 

Sadly, today the belief that the incompleteness Theorems even can
apply to \emph{FSOA} and, above all, that they have as a corollary
the semantic incompleteness of the \emph{full }second-order Logic
is almost unanimous. This can be seen either in \emph{Wikipedia} or
in the most specialized paper. Even in the introductory note of the
aforementioned document, Goldfarb writes:
\begin{quotation}
Finally, Gödel considers categoricity and syntactic completeness in
the setting of higher-order logics. {[}...{]} Noting then that Peano
Arithmetic is categorical - where by Peano Arithmetic he means the
second-order formulation - Gödel infers that if higher-order logic
is {[}semantically{]} complete, then there will be a syntactically
complete axiom system for Peano Arithmetic. At this point, he announces
his incompleteness theorem: \textquotedblleft The Peano axiom system,
with the logic of Principia mathematica added as superstructure, is
not syntactically complete\textquotedblright . He uses the result
to conclude that there is no (semantically) complete axiom system
for higher-order logic%
\footnote{W. Goldfarb, \emph{Note to {*}1930c}, \emph{op. cit}, vol. 3, pp.
14-15. %
}.
\end{quotation}
So interpreting, without a doubt, that Gödel refers to the second-order
categorical Arithmetic. But in this case the incompleteness Theorem
could not be applied! Goldfarb neglects that the \emph{full} second-order
induction principle, the only capable to ensures categoricity, generates
an uncountable quantity of axioms, so that the effective axiomatizability
of the system is not verified. As a matter of fact, the \emph{FSOA}
theory could be syntactically complete.

Sometimes, the semantic incompleteness of the \emph{full} second-order
Logic is ``concluded'' by an alternative approach to passing by
the (alleged) syntactic incompleteness of \emph{FSOA} (a \emph{Freudian}
stimulus?). It is supposed, by contradiction, that the valid statements
of the \emph{full} second-order are \emph{effectively enumerable}
theorems and then, by applying the incompleteness (or Tarski's) Theorem,
an absurdity is obtained%
\footnote{Emblematic examples, respectively, in: J. Hintikka, \emph{On Gödel},
Wadsworth Philosophers Series (2000), p. 22 and S. Shapiro, \emph{Do
Not Claim Too Much: Second-order Logic and First-order Logic}, Philosophia
Mathematica n. 3, vol. 7 (1999), p. 43.%
}. But this assumption is evidently excessive: semantic completeness
of a system, simply requires that all the valid statements are \emph{theorems},
not necessarily \emph{effectively enumerable} \emph{theorems}. In
the present case, conversely, we already know that this latter condition
is impossible, since, as it has been shown, the \emph{FSOA} is an
\emph{inherently semantic} theory. As a consequence, these proofs
really do not conclude the semantic incompleteness of the \emph{full}
second-order Logic (as, however, it can be done by the \emph{L-S}
Theorem or the invalidity of the compactness Theorem).

The main reasons of this unfortunate misunderstanding are due probably
to ambiguities of the used terminology, both ancient and modern.

\section{Clearing up the terms}

The expression ``extended predicate calculus'' is for the first
time used by Hilbert in the first edition (1928) of the aforementioned
\emph{Grundzüge der theoretischen Logik} where, with no doubt, indicates
the \emph{full} second-order Logic. This one was considered for the
first time in the \emph{Principia Mathematica} (\emph{PM}). The belief
that Gödel, in the aforementioned phrase, refers to \emph{FSOA} (explanation
\emph{b.}), implies that he with ``extended \emph{functional} calculus''
intends the same thing. But in which work he has shown or at least
suggested that the incompleteness Theorems can apply to the \emph{full}
second-order Logic? In none. 

In his proof of 1931, Godel refers to a formal system with a language
that, in addition to the first-order classical logic, allows the use
of non-bound \emph{functional variables} (ie without the possibility
to quantify on them)\emph{}%
\footnote{K. Gödel,\emph{ op. cit.}, vol. 1, p. 187. The purpose is simply to
be able to express formulas such as $\forall x(f(x))$, where \emph{f}
is a recursive function, which strictly are not permitted in the first-order
classical Logic.%
}. Then he proves that this is not a real extension of the language,
able, in particular, to hinder the applicability of the semantic completeness
Theorem. In the \emph{1932b} publication, Gödel declares the validity
of the incompleteness Theorems for a formal system (\emph{Z}), based
on first-order logic, with the axioms of Peano and an induction principle
defined by a recursive function. Certainly not a \emph{full} induction.
He adds: 
\begin{quotation}
If we imagine that the system \emph{Z} is successively enlarged by
the introduction of variables for classes of numbers, classes of classes
of numbers, and so forth, together with the corresponding comprehension
axioms, we obtain a sequence (continuable into the transfinite) of
formal systems that satisfy the assumptions mentioned above {[}...{]}%
\footnote{K. Gödel,\emph{ op. cit.}, vol. 1, p. 237.%
}
\end{quotation}
speaking explicitly of \emph{comprehension axioms} and \emph{formal}
systems. Finally, in the publication of 1934, which contains the last
and definitive proof of the first incompleteness Theorem, Gödel, having
the aim both to generalize and to simplify the proof, allows the quantification
either on the functional or propositional variables: a declared type
of second-order. However, appropriate comprehension axioms limit at
infinite countable the number of propositions%
\footnote{K. Gödel, \textit{On undecidable propositions of formal mathematical
systems}\textit{\emph{,}} \emph{op. cit.}, vol. 1, p. 353-354.%
}. Gödel never misses an opportunity to point out carefully that always
is referring to a formal system and that, consequently, the statements
are enumerable:
\begin{quotation}
Different formal systems are determined according to how many of these
types of variables are used. We shall restrict ourselves to the first
two types; that is, we shall use variables of the three sorts \emph{p},
\emph{q},\emph{ r},...{[}\emph{propositional variables}{]}; \emph{x},
\emph{y}, \emph{z},...{[}\emph{natural numbers variables}{]}; \emph{f},
\emph{g}, \emph{h}, ...{[}\emph{functional variables}{]}. We assume
that a denumerably infinite number of each are included among the
undefined terms (as may be secured, for example, by the use of letters
with numerical subscripts).

{[}...{]}

For undefined terms (hence the formulas and proofs) are countable,
and hence a representation of the system by a system of positive integers
can be constructed, as we shall now do%
\footnote{\emph{Op. cit.}, p. 350 and p. 355.%
}.
\end{quotation}
Therefore, we are in case 1) of the third section: far from \emph{full
}second-order. Yet, in the introduction to the same paper, Kleene,
in summarizing the work of Gödel, does not avoid commenting ambiguously:
\begin{quotation}
Quantified propositional variables are eliminable in favor of function
quantifiers. Thus the whole system is a form of full second-order
arithmetic (now frequently called the system of ``analysis'')%
\footnote{S. C. Kleene, \emph{op. cit.}, p. 339.%
}.
\end{quotation}
But he can mean only that the whole system is a \emph{formal version
}(perhaps as large as possible) of the \emph{full} second-order Arithmetic.
Maybe is exactly this the ``extended functional calculus'' to which
Gödel was referring in the examined words at the Congress? We will
discuss it in the next section.

Another source of mistake is probably related to use of the term \emph{metamathematics}.
Although Godel intends it in the modern broad sense that includes
any kind of argument beyond to the coded formal language of Mathematics
(and so with the possibility of using inherently semantic inferences
and/or making use of the concept of truth), in his theorems he employs
this term always limiting it to a \emph{formalizable} (though often
not yet formalized) use deductive and, indeed, even decidable. In
the short paper that anticipates his incompleteness Theorems%
\footnote{K. Gödel, \textit{op. cit.,} doc. \emph{1930b}, vol. 1, p. 143.%
}, for example, Gödel invokes a metamathematics able to decide whether
a formula is an axiom or not :
\begin{quotation}
{[}...{]}

IV. Theorem I {[}\emph{first incompleteness Theorem}{]} still holds
for all $\omega$-consistent extensions of the system \emph{S} that
are obtained by the addition of infinitely many axioms, provided the
added class of axioms is decidable, that is, provided for every formula
it is metamathematically decidable whether it is an axiom or not (here
again we suppose that in metamathematics we have at our disposal the
logical devices of \emph{PM}).

Theorems I, III {[}\emph{as the IV, but the added axioms are finite}{]},
and IV can be extended also to other formal systems, for example,
to the Zermelo-Fraenkel axiom system of set theory, provided the systems
in question are $\omega$-consistent.
\end{quotation}
But in both the subsequent rigorous proofs, he \emph{will formalize}
this process, which now is called metamathematical, using the recursive
functions, so revealing that, in the words just quoted, he refers
to the usual ``mechanical'' decidability. By the same token, even
in the theorem that concludes the consistency of the axiom of choice
and of the continuum hypothesis with the other axioms of the formal
Set Theory, he does the same: he uses the metamathematics only as
a simplification, stating explicitly that all ``the proofs could
be formalized'' and that ``the general metamathematical considerations
could be left out entirely''%
\footnote{K. Gödel, \textit{op. cit.,} doc. \emph{1940}, vol. 2, p. 34. %
}.

\section{An alternative explanation}

As noted, Gödel has never put in writing that his proofs of incompleteness
may be applied to the uncountable \emph{full} second-order Arithmetic
and furthermore it looks absolutely not reasonable to believe that
he deems it%
\footnote{I myself have changed my opinion reported in the book {[}12{]} and
in the paper {[}13{]} after a deeper analysis of the \emph{Kurt Gödel\textquoteright s
collected works}.%
}. In this section, therefore, we will examine the other possibility,
namely the \emph{a.} of the fourth section. It, remember, pretends
that Godel in 1930 believed, mistakenly, categorical a kind of formal
arithmetic and, in consequence of his incompleteness Theorems, semantically
incomplete its language. Is this reasonable (or more reasonable than
the previous case)?

Certainly not for the system considered by Gödel in his first proof
of 1931: in fact, the semantic completeness Theorem applies to it,
as Gödel himself remarks in note no. 55 of the publication%
\footnote{\emph{Op. cit.}, p. 187.%
}. Indeed, this is the first time in which the existence of non-standard
models for a formal arithmetic is proved: why Gödel does not report
it? The topic deserves a brief analysis.

As we have observed in third section, apart from the use of the incompleteness
Theorems, the existence of non-standard models for \emph{PA} can be
proved by the compactness Theorem, the \emph{L-S} \emph{upward} one,
or a theorem proved by Skolem in 1933 (cited in note no. 25); actually,
this last one resolves that \emph{PAT} is non-categorical, what \emph{a
fortiori} applies to \emph{PA}. The compactness Theorem is due precisely
to Gödel (1930) and descends from his semantic completeness Theorem;
but in none of his works Gödel ever uses it%
\footnote{This is confirmed by Feferman\emph{, op. cit.}, vol. 1, note n. 23,
p. 33.%
}. Moreover, despite its fundamental importance for the \emph{model
theory}, nobody - except Maltsev in 1936 and 1941 - uses it before
1945%
\footnote{This is attested both by Vaught and Fenstad: \emph{op. cit.}, vol.
1, p. 377 and vol. 2, p. 309.%
}. Not much more fortunate is the story of the \emph{L-S} Theorem\emph{.}
The first proof, by Löwenheim (1915), will be simplified by Skolem
in 1920%
\footnote{The respective works, relative to the first-order classical Logic
instead of the semantic completeness, are in: J. V. Heijenoort, \emph{From
Frege to Gödel}, Haward University Press (1967), pp. 232 - 251 and
pp. 254 - 263.%
}. In both cases, these theorems are \emph{downward} versions, able
to conclude the non-categoricity of the formal theory of the real
numbers and of the formal Set Theory, but not of \emph{PA}. However,
Skolem and Von Neumann suspect a much more general validity of the
result%
\footnote{Here is an especially interesting remark by Von Neumann in 1925, from
\emph{An axiomatization of set theory}, in \emph{From Frege to Gödel},
\emph{op. cit.}, p.\emph{ }412:
\begin{quotation}
{[}...{]} no categorical axiomatization of set theory seems to exist
at all; {[}...{]} And since there is no axiom system for mathematics,
geometry, and so forth that does not presuppose set theory, there
probably cannot be any categorically axiomatized infinite systems
at all.\end{quotation}
}. It seems that also Tarski was interested to this argument at that
time, probably getting the \emph{upward} version of the Theorem in
a seminar of 1928%
\footnote{This information is due to Maltsev: in his publication of 1936 (cited
in a next note), he claims to have known it by Skolem.%
}. In any case, the argument continues to have low popularity%
\footnote{It is significant, for example, that Hilbert did not mention this
theme in the seminary of Hamburg in 1927: \emph{The foundations of
mathematics}, in \emph{From Frege to Gödel}, \emph{op. cit.}, p.\emph{
}464. %
}, at least until the generalization of Maltsev in 1936%
\footnote{A. Maltsev, \emph{Untersuchungen aus dem Gebiete...}, Matematicheskii
sbornik 1, 323-336; English translation in: \emph{The metamathematics
of algebraic systems: collected papers} \emph{1936-1967}, Amsterdam
(1971). %
}, which, including for the first time the \emph{upward} version, will
allow the general conclusion that all the theories equipped with an
infinite model and a semantically complete language are not categorical.

In this context of disinterest for the topic, Gödel not only is no
exception, but his notorious Platonist inclination pushes him to distrust
and/or despise any interpretation that refers to objects foreign to
those that he believes existing independently of the considered theory.
Which, in all plausibility, also believes unique. As a matter of fact,
in the introduction of his first paper on the semantic completeness,
he shows to believe categorical even the first-order formal theory
of the real numbers%
\footnote{K. Gödel, \textit{op. cit.,} doc. \emph{1929}, vol. 1, pp. 61-63.%
}. 

On this basis, one can surmise the following alternative for the option
\emph{a.} When he discovers the non-categoricity of the formal arithmetical
system where his original incompleteness Theorems are applied, Gödel
is not so glad and immediately looks for an extension that, though
formal, is able to ensure the categoricity. Probably he believes to
have identified it in a \emph{formal version} of the \emph{full} second
order Arithmetic: just that one that will be considered in his generalized
proof of the first incompleteness Theorem of 1934%
\footnote{\emph{K. Gödel, op. cit.}, doc. \emph{1934}, vol. 1, p. 346.%
}, where quantification on the functional and propositional variables
are allowed, while respecting the formality. This hypothesis is consistent
with the fact that Gödel could admit the possibility that this language
is semantically incomplete, since in both versions of his semantic
completeness Theorem, he does not allow the use of quantifiers on
functional variables%
\footnote{\emph{Op. cit.}, doc.\emph{ 1929}, \emph{1930} and \emph{1930a}, vol.
1, p. 69, p. 121 and p. 125.%
}. Having in mind this generalization (perfectly definable as ``extended
functional calculus''), it would explain why in the meantime he communicates
the general result without mentioning the discovery of non-standard
models. For example, just after his famous announcement at the Congress
of Königsberg on September 6, 1930, Gödel declares%
\footnote{\emph{Op. cit.}, doc. \emph{1931a}, vol. 1, p. 203.%
}:
\begin{quotation}
(Assuming the consistency of classical mathematics) one can even give
examples of propositions (and in fact of those of the type of Goldbach
or Fermat) that, while contentually true, are unprovable in the formal
system of classical mathematics. Therefore, if one adjoins the negation
of such a proposition to the axioms of classical mathematics, one
obtains a consistent system in which a contentually false proposition
is provable.
\end{quotation}
Admitted the soundness, here he evades to mention standard models,
perhaps because he plans a more general proof than that one currently
available, valid for a formal arithmetic believed categorical. A theory
that, as announced in the main affirmation under scrutiny, uses the
``extended functional calculus'': inevitably, a semantically incomplete
language due to the syntactic incompleteness together with the (alleged)
categoricity. 

Anyway, the Skolem's proof of 1933 should have amazed him: not even
the system, semantically and syntactically complete, of the propositions
true for the standard model, is categorical. A first disturbing evidence
that the non-categoricity covers all formal (non-trivial) systems,
regardless of the syntactic completeness or incompleteness. Gödel,
in reviewing the Skolem's paper, only observes laconically - finally!
- that a consequence of this result, that is the non-categoricity
of \emph{PA}, was already derivable from his incompleteness Theorems%
\footnote{\emph{Op. cit.}, doc. \emph{1934c}, vol. 1, p. 379.%
}. Later, in any work (not only in the cited generalization of 1934),
he always will ignore the issue of the categoricity\emph{.} Equally,
\emph{apart from the questioned phrase of 1930, he never will return
to state that by his incompleteness Theorems can be derived the semantic
incompleteness of some language}.

Finally, the Henkin's Theorem of 1949 will prove that in every formal
system (and so, anywhere the incompleteness Theorems could be applied)
there is semantic completeness of the language and therefore, if infinite
models exist, there can not be categoricity.

\section{Conclusions}

We summarize briefly the conclusions of this paper:
\begin{enumerate}
\item The Richard's paradox can be interpreted as a meta-proof that the
semantic definitions are not countable, ie that they are conceivably
able to define each element of a set with cardinality greater than
the enumerable one (so, each real number). Further, the Berry's paradox
shows that a finite amount of symbols, differently interpreted, is
able to define an infinite amount of objects.
\item An axiomatic system that allows us to express and study an uncountable
amount of properties cannot be formal.
\item The categorical model of the \emph{full} second-order Arithmetic is
an interpretation that, before satisfying the axioms, associates 2\textsuperscript{$\aleph_{0}$}
different propositions to at least one sentence. The result, thus,
is a \emph{non-formal} axiomatic system.
\item An axiomatic system can be semantically complete despite employing
a semantically incomplete \emph{language}. In particular, the \emph{full}
second-order Arithmetic could be syntactically (and therefore also
syntactically, being categorical) complete.
\item The grouping of the axiomatic theories based on the expression order
is, in general, misleading their fundamental logical properties: these
are only a consequence of the premises. The essential tool for classification
remains the accomplishment of the Hilbertian formality.
\item The text of the Gödel's communication at the conference in Königsberg
on September 6, 1930 (never published by him), contains a mistake.
In the common understanding, not only this error is not reported,
but also is wrongly deduced, by it, that: a) the incompleteness Theorems
can also be applied to the categorical Arithmetic founded on the \emph{full}
second-order Logic; b) the semantic incompleteness of the \emph{full}
second-order Logic is a consequence of the incompleteness Theorems
(or of the ensuing Tarski's truth-undefinability Theorem).
\item On the basis of the publications of Gödel and of logic, the previous
interpretation is untenable.
\item As an alternative interpretation of the manuscript in question, it
is possible that Gödel is referring to the formal arithmetic considered
in his proof of 1934, in which the quantification on the functional
and propositional variables is allowed. If so, in 1930 he believed
categorical this theory and, as a consequence of its syntactic incompleteness,
supplied with a semantically incomplete language. This explanation
is consistent with the fact that both his original semantic completeness
Theorems cannot be applied to this system, due to the said quantification.
This view also explains why, being with the years more and more evident
the difficulty for the condition of categoricity, he never will repeat
similar affirmations.\\
On the other hand, he never corrected the phrase, plausibly because
he might not have worry about correcting an unpublished text. \end{enumerate}


\begin{thebibliography}{10}
\bibitem{key-4}ENDERTON H. B., \emph{Second-order and Higher-order
Logic}, Standford Encyclopedia of Philosophy (2007), http://plato.stanford.edu/entries/logic-higher-order/

\bibitem{key-23}FEFERMAN S. ET AL.\emph{, Kurt Gödel collected works},
ed. Feferman \emph{et al.}, Oxford univ. press (1986 -1995).

\bibitem{key-22}FENSTAD J. E., \emph{Selected works in logic}, ed.
Jens E. Fenstad, Univ. of Oslo (1970).

\bibitem{key-21}FREGE, \emph{Alle origini della nuova logica}, ed.
Boringhieri Torino (1983).

\bibitem{key-20}HEIJENOORT. V., \emph{From Frege to Gödel}, Haward
University Press (1967).

\bibitem{key-19}HENKIN L., \emph{The completeness of the first-order
functional calculus}, The journal of symbolic logic n. 14 (1949).

\bibitem{key-18}HILBERT D., ACKERMANN W., \emph{Principles of Mathematical
Logic}, ed. L. M. Hammond \emph{et al.}, Chelsea, New York (1950).

\bibitem{key-5}HINTIKKA J., \emph{On Gödel}, Wadsworth Philosophers
Series (2000).

\bibitem{key-6}KENNEDY J., \emph{Completeness Theorems and the Separation
of the First and Higher-Order Logic} (2008), http://igiturarchive.
library.uu.nl/lg/2008-0317-201019/UUindex.html 

\bibitem{key-7}LEWIS C. I., \emph{A Survey of Formal Logic}, Berkeley
University of California Press (1918).

\bibitem{key-8}MALTSEV A., \emph{The metamathematics of algebraic
systems: collected papers} \emph{1936-1967},\emph{ }ed. Benjamin Franklin
Wells, Amsterdam (1971).

\bibitem{key-9}RAGUNÌ G., \emph{I confini logici della matematica},
ed. Aracne, Roma (2010).

\bibitem{key-11}RAGUNÌ G., \emph{The Gödel\textquoteright s legacy:
revisiting the Logic}, Intellectual Archive n. 2, vol. 1 (2012).

\bibitem{key-24}RAGUNÌ G., \emph{La logica matematica dopo Gödel},
ed. Amazon (2013).

\bibitem{key-12}ROSSBERG M., \emph{First-Order Logic, Second- Order
Logic and Completeness}, ed. Hendricks \emph{et al}., Logos Verlag
Berlin (2004), http://www.standrews. ac.uk/\textasciitilde{}mr30/papers/RossbergCompleteness.pdf 

\bibitem{key-13}RUSSELL B., WHITEHEAD A. N., \emph{Principia Mathematica},
Cambridge University Press (1913).

\bibitem{key-14}SHAPIRO S., \emph{Do Not Claim Too Much: Second-order
Logic and First-order Logic}, Philosophia Mathematica n. 3, vol. 7
(1999).

\bibitem{key-15}TURING A. M., \emph{On computable numbers, with an
application to the Entscheidungsproblem}, Proceed. of the London Mathematical
Society n. 42, vol. 2 (1937).

\bibitem{key-16}WRIGHT C., \emph{On Quantifying into Predicate Position:
Steps towards a New(tralist) Perspective} (2007), http://philpapers.org/autosense.pl?searchStr=Crispin\%20Wright.

\bibitem{key-17}ZALAMEA F., \emph{Filosofía sintética de las Matemáticas
contemporáneas} (2009), http://files.acervopeirceano.webnode.es/200000065-
18c1b19bb9/Zalamea-Fil-Sint-Mat-Cont.pdf 13\end{thebibliography}
\end{document}